\documentstyle[11pt]{article}

\newcommand{\Bbb}{\bf}
\def\frak#1{{#1}}

\def\ket#1{| #1\rangle}		

\def\be{\begin{equation}}
\def\en{\end{equation}}
\def\bea{\begin{eqnarray}}
\def\ena{\end{eqnarray}}
\def\bean{\begin{eqnarray*}}
\def\enan{\end{eqnarray*}}

\newcommand{\gsl}{{\frak {sl}}}
\newcommand{\ggl}{{\frak {gl}}}
\newcommand{\slt}{\gsl_2}
\newcommand{\sln}{\gsl_n}
\newcommand{\glt}{\ggl_2}
\newcommand{\dyslt}{{{\cal D}{ Y}(\gsl_2)_k}}

\newcommand{\slh}{\widehat{\gsl}}

\newcommand{\del}{\partial}

\def\slth{\slh_2}

\def\ket#1{|#1\rangle}
\def\del#1#2{{}_{#1}\partial_{#2}}

\def\Z{{\Bbb Z}}

\def\C{{\Bbb C}}

\begin{document}

\begin{flushright}
YITP-96-10 \\
May 1996 
\end{flushright}
\vspace{24pt}

\begin{center}
\begin{LARGE}
{  Free Field Representation of Level-$k$ \break
 Yangian Double 
$\dyslt$}\par
\vskip 1mm
{ and Deformation of Wakimoto Modules${}^{\#}$}
\end{LARGE}

\vspace{30pt}

\begin{large}
Hitoshi Konno\raisebox{2mm}{$\dagger$} 
\end{large}

\vspace{6pt}

\begin{large}
{\it Yukawa Institute for Theoretical Physics,} \\
{\it Kyoto University, Kyoto 606-01, Japan\raisebox{2mm}{$\star$}}
\end{large}

\vspace{35pt}

{ABSTRACT}
\end{center}

\vspace{20pt}
Free field representation of level-$k (\not=0,-2)$ Yangian 
double $\dyslt$ and a corresponding 
deformation of Wakimoto modules are
presented. We also realize two types of vertex operators 
intertwining these modules.

\vspace{12pt}
${\#}$ Accepted for publication in Lett.Math.Phys.
\vspace{10pt}
\vfill
\hrule

\vskip 3mm
\begin{small}

\noindent\raisebox{2mm}{$\dagger$} Yukawa fellow

\noindent\raisebox{2mm}{$\star$} E-mail: konno@yukawa.kyoto-u.ac.jp

\end{small}

\setcounter{section}{0}
\setcounter{equation}{0}
\section{Introduction}
Yangian is a Hopf algebra associated with the rational solutions of
the quantum Yang-Baxter equation(YBE)\cite{Dr1}. 
In \cite{LS,Smir,BL}, it has been shown that the Yangian
double ${\cal D}Y(g)$, a quantum double of the Yangian\cite{Dr2}, 
is a relevant
object to describe a mathematical structure of massive integrable 
quantum field theory. However, this idea has not yet been realized
successfully. It seems that this is because only level-0 
Yangian double has been considered, 
where infinite dimensional representations are lacking.
This situation makes a contrast to the success in integrable spin chains,
where the quantum affine algebras,
a Hopf algebra associated with the trigonometric solutions of
the quantum YBE, have been applied\cite{FR,JM}. 

Recently, a breakthrough has been brought by Iohara, Kohno\cite{IK,Ioh}
and Khoroshkin, Lebedev, Pakuliak\cite{Kho,KLP}. 
They have succeeded to define a central
extension of the Yangian double ${{\cal D}Y(\sln)}$ 
by means of the Drinfeld currents\cite{Dr3} and to
obtain infinite dimensional representations which are analogous to 
the level-1 highest weight representations of the affine Lie algebra
$\widehat{\sln}$ and the quantum affine algebra $U_{q}(\widehat{sl_n})$.

In this letter, we extend this line to the higher level representation
for  the central extension of ${{\cal D}Y(\slt)}$.
We consider a representation of the level-$k\ (\not=0,-2)$ Drinfeld
currents in terms of the three free bosonic fields $\Phi, \phi$ and
$\chi$. This can be regarded as a deformation of the so-called
Wakimoto construction of the affine Lie algebras\cite{Wak}.
As a consequence, an analogue of the level-$k$ Wakimoto modules are
obtained as a certain restriction of the Fock spaces of these bosons. 
We also give a boson representation of vertex operators intertwining
the resultant modules. 

This paper is organized as follows. In section 2,
we briefly recall the definition of the central extension of the 
 Yangian double ${{\cal D}Y(\slt)}$.
In section 3, introducing three bosonic fields we give a free field
construction of the level-$k$ Drinfeld currents.
In section 4, we point out the existence of a couple of operators 
called screening charges which commutes with all the currents. 
Their free field representation is also given.
In section 5, we consider vertex operators called type
I and type II in the terminology by Jimbo and Miwa\cite{JM} and give their
bosonized expressions. The final section is devoted to discussions.
 
\setcounter{section}{1}
\setcounter{equation}{0}

\section{Central extension of the Yangian double}


\noindent
{\bf Definition 2.1}
The central extension of the yangian double ${{\cal D}Y(\slt)}$
\cite{IK,Kho} is a Hopf algebra over $\C[[\hbar]]$ generated by 
the derivation operator $d$, central element $c$ and the symbols $
e_m, f_m, h_m, \ m\in\Z$ satisfying the following relations.
\bea
&&[d, e(u)]=\frac{d}{du}e(u),\qquad 
[d, f(u)]=\frac{d}{du}f(u), \qquad
[d, h^{\pm}(u)]=\frac{d}{du}h^{\pm}(u),\nonumber \\
&& e(u)e(v)=\frac{u-v+\hbar}{u-v-\hbar}e(v)e(u),\nonumber \\
&& f(u)f(v)=\frac{u-v-\hbar}{u-v+\hbar}f(v)f(u),\nonumber \\
&& h^{\pm}(u)e(v)=\frac{u-v+\hbar}{u-v-\hbar}e(v)h^{\pm}(u),\label{relations} \\
&& h^+(u)f(v)=\frac{u-v-(1+c)\hbar}{u-v+(1-c)\hbar}f(v)h^+(u),\nonumber \\
&& h^-(u)f(v)=\frac{u-v-\hbar}{u-v+\hbar}f(v)h^-(u),\nonumber \\
&& [ h^{\pm}(u), h^{\pm}(v) ]=0,\nonumber \\
&& h^+(u)h^-(v)=\frac{u-v+\hbar}{u-v-\hbar}
        \frac{u-v-(1+c)\hbar}{u-v+(1-c)\hbar}h^-(v)h^+(u),\nonumber \\
&& [e(u), f(v)]=\frac{1}{\hbar}(\delta(u-(v+\hbar c))h^+(u)-\delta(u-v)h^-(v)),\nonumber
\ena
where  
\bea 
&& e^{\pm}(u)=\pm\sum_{m\geq 0\atop m<0}e_m u^{-m-1},\qquad
 f^{\pm}(u)=\pm\sum_{m\geq 0\atop m<0}f_m u^{-m-1},\qquad
\nonumber \\
&&  h^{\pm}(u)=1\pm\hbar\sum_{m\geq 0\atop m<0}h_m u^{-m-1}, 
\label{modes}
\ena
and 
\bea
&& e(u)=e^+(u)-e^-(u),\qquad f(u)=f^+(u)-f^-(u),\nonumber\\ 
&& \delta(u-v) = \sum_{n+m=-1}u^nv^m.\nonumber
\ena
The coproduct is given by\cite{Iohara} 
\bea
&& \Delta(d)=d\otimes 1 +1\otimes d,
\quad  \Delta(c)=c\otimes 1 +1\otimes c,\nonumber \\
&& \Delta(e^{\epsilon}(u))=e^{\epsilon}(u)\otimes 1\nonumber \\
&&\qquad+\Bigl[1\otimes 1+
\hbar^2 f^{\epsilon}(u+\hbar-\delta_{\epsilon,+}c \hbar)
\otimes e^{\epsilon}(u-\frac{1}{2}\epsilon c_1\hbar )\Bigr]^{-1}
h^{\epsilon}(u)\otimes e^{\epsilon}(u-\frac{1}{2}\epsilon c_1\hbar )
,\nonumber \\
&& \Delta(h^{\epsilon}(u))=
\Bigl[1\otimes 1+
\hbar^2 f^{\epsilon}(u+\hbar-\delta_{\epsilon,+}c \hbar)
\otimes e^{\epsilon}(u-\hbar-\frac{1}{2}\epsilon c_1\hbar )\Bigr]^{-2}
\nonumber\\
&&\qquad\qquad\times 
h^{\epsilon}(u)\otimes h^{\epsilon}(u-\frac{1}{2}\epsilon c_1\hbar ),
\label{comulti} \\
&& \Delta(f^{\epsilon}(u))=1\otimes
f^{\epsilon}(u+\frac{1}{2}c_1\hbar)\nonumber \\
&&\qquad+
\Bigl[1\otimes 1+
\hbar^2 f^{\epsilon}(u+\delta_{\epsilon,+}c_2 \hbar)
\otimes e^{\epsilon}(u-\hbar+\frac{1}{2} c_1\hbar
+\delta_{\epsilon,+}c \hbar )\Bigr]^{-1}\nonumber \\
&&\qquad\qquad\times f^{\epsilon}(u+\delta_{\epsilon,+}c_2 \hbar)\otimes 
h^{\epsilon}(u+\frac{1}{2}c_1\hbar+\delta_{\epsilon,+}c\hbar ),
\nonumber\ena
where 
$c_1=c\otimes 1,\ c_2=1\otimes c$, 
$\epsilon = \pm, \delta_{+,+}=1$ and $\delta_{+,-}=0$.

\vspace{5mm}\noindent
{\it Remark.} The formulae (\ref{relations}) and (\ref{comulti}) are 
converted to those in \cite{IK} by the following rule\cite{Iohara}.
\bea
&&H^{\pm}(u\mp \frac{1}{4}c\hbar)=h^{\pm}(u),\quad 
E(u)=e(u),\quad F(u+\frac{1}{2}c\hbar)=f(u).
\ena
Note also 
\bea
&&h^{\epsilon}(u)e^{\epsilon}(u+\hbar)=e^{\epsilon}(u-\hbar)h^{\epsilon}(u),
\nonumber \\
&&h^{\epsilon}(u)f^{\epsilon}(u-\hbar-\delta_{\epsilon,+}c\hbar)
=f^{\epsilon}(u+\hbar-\delta_{\epsilon,+}c\hbar)
h^{\epsilon}(u).
\ena

\setcounter{section}{2}
\setcounter{equation}{0}

\section{Free bosons}
Let  us introduce a Heisenberg algebra generated by 
$a_{\Phi,n},\ a_{\phi,n },\ a_{\chi,n}\ n\in\Z$ and $a_{\Phi}, 
a_{\phi}, \ a_{\chi}, \partial_{\Phi}, \partial_{\phi},  \partial_{\chi}$ 
satisfying the commutation relations
\bea
&& [a_{\Phi,m}, a_{\Phi,n} ]=\frac{k+2}{2}m\delta_{m+n,0},\qquad 
 [\partial_{\Phi}, a_\Phi]=\frac{k+2}{2}, \nonumber\\
&& [a_{\phi,m}, a_{\phi,n} ]=-m\delta_{m+n,0},\qquad  
[\partial_{\phi}, a_\phi]=-1, \label{heisenberg}\\
&& [a_{\chi,m}, a_{\chi,n} ]=m\delta_{m+n,0},\qquad  
[\partial_{\chi}, a_\chi]=1\nonumber\\
&&[a_{X,m},a_{X',n}]=0,   \qquad X\not=X'\nonumber
\ena
with $k\in \C(\not=0,-2)$.

Let $\ket{l;s,t}$ be a vacuum state with $\Phi, \phi, \chi$-charges
$l, -s, t$: 
\bea
&&\ket{l;s,t}=e^{\frac{l}{k+2}a_{{\Phi}}+sa_{\phi}+ta_{\chi}}\ket{0},\\
&& a_{X,n}\ket{0}=0,\quad n{>0}, \quad \partial_{X}\ket{0}=0\nonumber
\ena
for $X=\Phi,\phi,\chi$.
 Let ${\cal F}_{l,s,t}$ be a  Fock space constructed  on 
$\ket{l;s,t}$:
\bea
&&{\cal F}_{l,s,t}=\Bigl\{\prod a_{\Phi,-n}\prod a_{\phi,-n'}\prod a_{\chi,-n''}\ket{l;s,t}\Bigr\}
\ena
with $n, n', n''\in \Z_{>0}$. 

It is convenient to introduce  generating functions  $X( u; A,B)$,
$X=\Phi, \phi, \chi$
\bea
&& X(u; A,B) =\sum_{n>0}\frac{\ a_{X,-n}}{ n}(u+A\hbar)^{n}
              -\sum_{n>0}\frac{a_{X,n}}{n}(u+B\hbar)^{-n} \nonumber \\ 
&&	\qquad\qquad\qquad        +  \log(u+B\hbar)\partial_{{X}}+a_X.
\ena
We also use the notation $X(u; A,A)=X(u;A)$.
We denote by $:\ :$  the normal ordered product
\bea
&&:\exp\{X(u;A,B)\}:=\exp\Bigl\{
\sum_{n>0}\frac{a_{X,-n}}{n}(u+A\hbar)^n\Bigr\}e^{a_{X}}\nonumber\\
&&\qquad\qquad\times
(u+B\hbar)^{\partial_{X}}\exp\Bigl\{-\sum_{n>0}\frac{a_{X,n}}{n}(u+B\hbar)^{-n}\Bigr\}
\nonumber\ena
and by ${}_\alpha\partial_u$ the difference operator
\bea
&&\del{\alpha}{u}g(u)=\frac{g(u+\alpha\hbar)-g(u)}{\hbar}.
\nonumber\ena

\vspace{5mm}
\noindent
{\bf Definition 3.1}
 The currents $e(u), f(u)$ and $h^{\pm}(u)$ acting on 
${\cal F}_{l,s,t}$ are defined by 
\bea
&& e(u)=-:[\del{1}{u} \exp\{-\chi(u;-(k+2))\}]\exp\{-\phi(u; -(k+1),-(k+2)) \}:
\nonumber \\
\label{dre}\\
&& f(u)=\frac{1}{\hbar}
:\Bigl[
\exp\Bigl\{ \sum_{n>0}\frac{ a_{\Phi, n}}{n}[(u-2\hbar)^{-n}-
u^{-n}]\Bigr\}
\Bigl(\frac{u}{u-2\hbar}\Bigr)^{\partial_{\Phi}}\nonumber\\
&&\qquad\qquad \times\exp\{\phi(u;-1,0)+\chi(u;-1)\}\nonumber \\
&&\qquad\qquad -\exp\Bigl\{  \sum_{n>0}\frac{2a_{\Phi, -n}}{(k+2)n}
[(u-(k+3)\hbar)^{n}-
(u-\hbar)^{n}]\Bigr\}\nonumber \\
&&\qquad\qquad\times \exp\{\phi(u;-(k+3), -(k+2))+\chi(u; -(k+2))\}
\Bigr]:\label{drf}
\\
&& h^+(u)=
\exp\Bigl\{  \sum_{n>0}\frac{a_{\Phi, n}}{n}[(u-(k+2)\hbar)^{-n}-
(u-k\hbar)^{-n}]\Bigr\} 
\Bigl(\frac{u-k\hbar}{u-(k+2)\hbar}\Bigr)^{\partial_{\Phi}}
\nonumber \\
&& \qquad\qquad\times\exp\Bigl\{  \sum_{n>0}\frac{a_{\phi, n}}{n}
[(u-(k+2)\hbar)^{-n}-
(u-k\hbar)^{-n}]\Bigr\}
\Bigl(\frac{u-k\hbar}{u-(k+2)\hbar}\Bigr)^{\partial_{\phi}}
\nonumber \\
\label{drhp}\\
&& h^-(u)=\exp\Bigl\{  \sum_{n>0}\frac{2a_{\Phi, -n}}{(k+2)n}
[(u-(k+3)\hbar)^{n}-
(u-\hbar)^{n}]\Bigr\}
\nonumber \\
&& \qquad\qquad \times \exp\Bigl\{  \sum_{n>0}\frac{a_{\phi, -n}}{n}
[(u-(k+3)\hbar)^{n}-
(u-(k+1)\hbar)^{n}]\Bigr\}\label{drhm}
\ena
We also define the operator $d$ by
\bea
&& d=d_{\Phi}+d_{\phi}+d_{\chi},\label{derivation}\\
&& d_{\Phi}=\frac{2}{k+2}\Bigl(a_{\Phi,-1}\partial_{\Phi}+
\sum_{n\in\Z_{>0}}a_{\Phi,-(n+1)}a_{\Phi,n}\Bigr),\nonumber\\
&& d_{\phi}=-a_{\phi,-1}\partial_{\phi}-
\sum_{n\in\Z_{>0}}a_{\phi,-(n+1)}a_{\phi,n},
\nonumber\\
&&d_{\chi}=a_{\chi,-1}\partial_{\chi}+
\sum_{n\in\Z_{>0}}a_{\chi,-(n+1)}a_{\chi,n}.\nonumber
\ena

\noindent
{\bf Proposition 3.1} {\it For $X=\Phi,\phi, \chi$,}
\bea
&&e^{\gamma d}a_{X,-m}e^{-\gamma d}=\sum_{n\geq 0}
\frac{(m+n-1)!}{(m-1)!n!}\gamma^na_{X, -(m+n)},\qquad n\geq 1
\nonumber\\
&&e^{\gamma d}a_{X,m}e^{-\gamma d}=\sum_{0\leq n <m }
\frac{(-)^nm!}{(m-n)!n!}\gamma^na_{X, m-n}+(-)^m\gamma^m\partial_{X},
\qquad n\geq 1\nonumber\\
&&e^{\gamma d}a_{X}e^{-\gamma d}=a_{X}+\sum_{n\geq 1}
\frac{a_{X, -n}}{n}\gamma^n,\qquad 
e^{\gamma d}\partial_{X}e^{-\gamma d}
=\partial_{X}\nonumber
\ena

By direct calculation using (\ref{dre})$\sim$(\ref{derivation}), we have 

\vspace{5mm}
\noindent
{\bf Theorem 3.1}
{\it By analytic continuation, the operator $d$ and the 
currents $e(u), f(u)$ 
$h^{\pm}(u)$ satisfy the same relations as (\ref{relations}) with $c=k$ on 
${\cal F}_{l,s,t}$.}

\vspace{5mm}
We denote by $\dyslt$ the algebra generated by the currents $e(u), f(u), 
h^{\pm}(u)$ and $d$ on ${\cal F}_{l,s,t}$.
 
\vspace{5mm}\noindent
{\it  Remark.}\
In the classical limit $\hbar\to 0$, the currents $e(u), f(u),
h^{\pm}(u)$ 
tend to
the bosonized form of the level-$k$ 
currents for the affine Lie algebra $\widehat{\slt}$:
\bea
&& e(u)\to -\Bigl[\partial e^{-\chi(u)}\Bigr] e^{-\phi(u)},\nonumber\\
&& f(u)\to \Bigl[(k+2)\partial \phi(u)+(k+1)\partial\chi(u)+2\partial
\Phi(u)\Bigr]e^{\chi(u)+\phi(u)},\nonumber\\
&& \frac{1}{\hbar}\Bigl( h^+(u)-h^-(u)\Bigr) \to  2\Bigl(\partial \phi(u)+\partial\Phi(u)\Bigr)
\nonumber\ena
with
\bea
&& X(u)=a_{X}+(\log u)\partial_{X}-\sum_{n\in \Z(\not=0)}
\frac{a_{X,n}}{n}u^{-n}, \qquad X=\Phi,\phi,\chi.
\nonumber
\ena
In the same limit, the operator $d$ is nothing but the $L_{-1}$
operator of the Virasoro algebra constructed via Sugawara form from the above
classical currents. 

For the level-$k$ quantum affine
algebra $U_q(\widehat{sl_2})$,
an analogous free field representation was obtained 
in \cite{Shir,Matsuo1,ABG,KQS}.

\setcounter{section}{3}
\setcounter{equation}{0}

\section{Deformed Wakimoto modules}

It is easy to show

\vspace{5mm}
\noindent
{\bf Proposition 4.1}
\bea
&&[e(u), \partial_{\phi}+\partial_{\chi}]=0,\quad 
[h^{\pm}(u), \partial_{\phi}+\partial_{\chi}]=0,\quad 
[f(u), \partial_{\phi}+\partial_{\chi}]=0.\nonumber
\ena 
We hence restrict the Fock space ${\cal F}_{l,s,t}$ to the $s=t$ sector
without loss of generality.

\vspace{5mm}
\noindent
{\it Remark.}\  
On ${\cal F}_{l,s,s}$, the currents $e(u), f(u)$  and $h^{\pm}(u)$
are single valued so that the expansion such as (\ref{modes}) makes sense.

\vspace{5mm} 
Let us next consider the following operators.
\bea
&&\xi(u)=:\exp\{-\chi(u; -(k+2)) \}:, \qquad \eta(u)=:\exp\{\chi(u; -(k+2))\}: 
\nonumber
\ena
We have 
\bea
&&\xi(u)\eta(v)=-\eta(v)\xi(u)\sim\frac{1}{u-v}.
\label{xieta}
\ena
Here $\sim$ means that the relation holds modulo regular terms.
The fields $\xi(u)$ 
and $\eta(u)$ are single valued on ${\cal F}_{l,s,s}\ s\in\Z$. 
From (\ref{xieta}),
the zero-modes $\xi_0=\oint\frac{du}{2\pi i u}\xi(u)$ and 
$\eta_0=\oint\frac{du}{2\pi i }\eta(u)$ anticommute
$ \{\xi_0,\eta_0\}=0$. Note also $\xi_0^2=0=\eta_0^2$.
In addition,  the following equations hold in the sense of analytic continuation.

\vspace{5mm}
\noindent
{\bf Proposition 4.2}
\bea
&& h^{\pm}(u)\eta(v)=\eta(v)h^{\pm}(u)\sim 0,\nonumber\\
&& f(u)\eta(v)=-\eta(v)f(u)\sim 0,\nonumber\\
&& e(u)\eta(v)=-\eta(v)e(u)\nonumber \\
&&\qquad\qquad\sim \del{1}{v}
\Bigl(\frac{1}{u-v+\hbar}\exp\{-\phi(u;-(k+2),-(k+3))\}\Bigr).\nonumber
\ena
Therefore the zero-mode $\eta_0$ commutes
with the action of  $\dyslt$.

From this we restrict the Fock space ${\cal F}_{l,s,s}$ to the kernel of $\eta_0$. We hence arrive at the deformation of the Wakimoto modules
\bea
&& {\cal F}_l=\bigoplus_{s\in \Z}{\rm Ker}( \eta_0: {\cal
F}_{l,s,s}\to 
{\cal F}_{l,s,s+1}).
\ena

For the level$-k$ $U_q(\widehat{sl_2})$, the $q-$deformation of the
Wakimoto modules were obtained in \cite{Matsuo2,Kon1,Kon2}.

\vspace{5mm}
\noindent
{\it Remark.}\ For $m\geq 0$
\bea
e_{m}\ket{l;0,0}&=&\oint\frac{du}{2\pi i }u^m e(u)\ket{l;0,0}\nonumber\\
&=&\oint\frac{du}{2\pi i }u^m\frac{1}{\hbar}\Bigl(
e^{-\sum_{n>0}\frac{a_{\chi,-n}}{n}(u-(k+1)\hbar)^n}-
e^{-\sum_{n>0}\frac{a_{\chi,-n}}{n}(u -(k+2)\hbar)^n}
\Bigr)
\nonumber \\
&&\qquad\times 
e^{-\sum_{n>0}\frac{a_{\phi,-n}}{n}(u-(k+1)\hbar)^n}\ket{l;0,0}
\nonumber\\
&=&0.
\ena
In the same way, for $m\geq 0$,
\bea
&&f_{m}\ket{l;0,0}=\delta_{m,0}2l\ket{l;0,0}
\nonumber \\
&&\quad+\sum_{a=2}^{\infty}
\Bigl(\matrix{l+a-1\cr a\cr}\Bigr)2^a\hbar^{a-1}\nonumber\\
&&\qquad\times
\oint\frac{du}{2\pi i}u^{m-a}
e^{\sum_{n>0}\frac{a_{\chi,-n}}{n}(u-\hbar)^n
+\sum_{n>0}\frac{a_{\phi,-n}}{n}(u-\hbar)^n}\ket{l;0,0},
\ l\not=0, \nonumber \\
&&f_{m}\ket{0;0,0}=0,\nonumber\\
&&h_m\ket{l;0,0}=\delta_{m,0}\frac{2l}{k+2}\ket{l;0,0}
\nonumber \\
&&\qquad\qquad+
(1-\delta_{m,0})((k+2)\hbar)^m\sum_{a=0\atop m+1\geq a}^{l}
\frac{l(l+m-a)!}{a!(l-a)!(m-a+1)!}\Bigl(\frac{-k}{k+2}\Bigr)^a\ket{l;0,0}.
\nonumber
\ena
Hence the states $\ket{l;0,0}$ are not the highest weight states in the usual sense except for the case $l=0$.
However, $\ket{l;0,0}\ l\not=0$ recovers  the highest weight state conditions in the limit $\hbar\to 0$. We hence identify the state 
$\ket{l;0,0}$ with the highest weight state $\ket{\lambda_{l}}$ in the 
affine Lie algebra $\slth$,
where $\lambda_l=(k-l)\Lambda_0+l\Lambda_1$ are dominant integral
weights of level$-k$ with $\Lambda_0$ 
and $\Lambda_1$ being the fundamental weights.  This identification is
natural, because the zero-modes of the free bosons (\ref{heisenberg}) 
are not affected 
by the deformations.
 
Next we consider the screening operator defined by
\bea
&&S(u)_{[J]}=:[\del{1}{u}\exp\{-\chi(u;-1)\}]\exp\{-\phi(u;-1,0)\}:
\nonumber\\
&&\qquad\qquad\times
\exp\Bigl\{-\sum_{n>0}\frac{2a_{\Phi,-n}}{(k+2)n}
(u-\hbar)^n\Bigr\} e^{-\frac{2}{k+2}a_{\Phi}}
\prod_{j=0}^{J}
\Bigl(\frac{u-(2+(k+2)j)\hbar}{u-(k+2)j\hbar}\Bigr)^{\partial_{\Phi}}
\nonumber \\
&&\qquad\qquad\times
\exp\Bigl\{\sum_{j=0}^J\sum_{n>0}
\frac{a_{\Phi,n}}{n}[(u-(k+2)j\hbar)^{-n}-(u-(2+(k+2)j)\hbar)^{-n}\Bigr\}
\nonumber \\
\ena
with $J\in\Z_{>0}$.

We have

\vspace{5mm}
\noindent
{\bf Proposition 4.3}
\bea
&& h^{+}(u)S(v)_{[J]}=S(v)_{[J]}h^{+}(u)\sim 0,\nonumber\\
&& e(u)S(v)_{[J]}=S(v)_{[J]}e(u)\sim 0,\nonumber\\
&& S(u)_{[J]}\eta(v)=\eta(v)S(u)_{[J]}\nonumber\\
&& \sim -\del{1}{v}
\Bigl(\frac{1}{u-v+(k+2)\hbar}\exp\{-\phi(v;-(k+3),-(k+2))\}
\nonumber \\
&&\times
\exp\Bigl\{-\sum_{n>0}\frac{2a_{\Phi,-n}}{(k+2)n}v^n\Bigr\}
e^{-\frac{2}{k+2}a_{\Phi}}
\prod_{j=0}^{J}\Bigl(\frac{v-(1+(k+2)j)\hbar}
{v-(-1+(k+2)j)\hbar}\Bigr)^{\partial_{\Phi}}\nonumber \\
&&\times
\exp\Bigl\{\sum_{j=0}^J\sum_{n>0}\frac{a_{\Phi,n}}{n}[
(v-(-1+(k+2)j)\hbar)^{-n}-(v-(1+(k+2)j)\hbar)^{-n}]\Bigr\}
\Bigr),\nonumber\\
&&[S(u)_{[J]},\partial_{\phi}+\partial_{\chi}]=0.\nonumber
\ena
{\it In addition, in the limit $J\to \infty$,}
\bea
&& h^{-}(u)S^{}(v)=S^{}(v)h^{-}(u)\sim 0,\nonumber\\
&& f(u)S(v)=S(v)f(u)\nonumber \\
&&\qquad \sim -\del{k+2}{v}
\Bigl(\frac{1}{u-v}
\exp\Bigl\{-
\sum_{n>0}\frac{2a_{\Phi,-n}}{(k+2)n}
(v-\hbar)^n\Bigr\}\nonumber \\
&&\qquad\qquad \times e^{-\frac{2}{k+2}a_{\Phi}}
\prod_{j=1}^{\infty}
\Bigl(\frac{v-(2+(k+2)j)\hbar}{v-(k+2)j\hbar}\Bigr)^{\partial_{\Phi}}
\nonumber \\
&&\quad\qquad\times
\exp\Bigl\{\sum_{j=1}^\infty\sum_{n>0}
\frac{a_{\Phi,n}}{n}[(v-(k+2)j\hbar)^{-n}-(v-(2+(k+2)j)\hbar)^{-n}]\Bigr\}
\Bigr),\nonumber
\ena
{\it where $S(u)=\lim_{J\to \infty}S(u)_{[J]}$.}

\vspace{5mm}
\noindent
Therefore the screening charge $S=\oint\frac{du}{2\pi i}S(u)_{[J]}$ commutes
with all the currents in $\dyslt$ and $\eta_0$
in the limit $J\to\infty$.
The charge $S$ yields a linear map $S: {\cal F}_{l}\to {\cal F}_{l-2}$.

\setcounter{section}{4}
\setcounter{equation}{0}
\section{Vertex operators}
Let $V^{(l)}=\oplus_{m=0}^l\C[[\hbar]] w_m$ be the 
$l+1$-dim representation of $\slt$. Let   
$V^{(l)}_u=V^{(l)}\otimes \C[[\hbar]][[u^{-1},u]$ be the evaluation module 
 on which the action of $\dyslt$ is defined by the following 
relations\cite{CP}. 
\bea
&&e_n(w_m\otimes u^a)=m
\Bigl(u+\frac{l-2m+1}{2}\hbar\Bigr)^n(w_{m-1}\otimes u^a)
,\nonumber \\
&&f_n(w_m\otimes u^a)=(l-m)\Bigl(u+\frac{l-2m-1}{2}\hbar\Bigr)^n
(w_{m+1}\otimes u^a),\label{evaluation}\\
&&h_n(w_m\otimes u^a)=\Bigl[(m+1)(l-m)\Bigl(u+\frac{l-2m-1}{2}\hbar\Bigr)^n
\nonumber \\
&&\qquad\qquad \qquad
-m(l-m+1)\Bigl(u+\frac{l-2m+1}{2}\hbar\Bigr)^n\Bigr]
(w_{m}\otimes u^a), \nonumber\\
&&d (w_{m}\otimes u^a)=-a(w_{m}\otimes u^{a-1})\nonumber
\ena
for $a, n\in \Z$, where we set $w_{-1}=w_{l+1}=0$.

Let $\lambda_{l}$ be the dominant integral weights. 

\vspace{5mm}
\noindent
{\bf Definition 5.1}
The type I 
and the type II vertex operators are the following intertwiners.
\bea  
{\rm Type\ I}&&\Phi_{\lambda_{l_1}}^{\lambda_{l_3} V^{(l)}}(u):
{\cal F}_{l_1}\to {\cal F}_{l_3}\otimes V^{(l)}_u,\nonumber\\
{\rm Type\ II }&&\Psi_{\lambda_{l_1}}^{\lambda_{l_3} V^{(l)}}(u):
{\cal F}_{l_1}\to V^{(l)}_u\otimes {\cal F}_{l_3} .\nonumber
\ena
satisfying
\bea
&&\Phi_{\lambda_{l_1}}^{\lambda_{l_3} V^{(l)}}(u)x=\Delta(x)
\Phi_{\lambda_{l_1}}^{\lambda_{l_3} V^{(l)}}(u),\nonumber \\
&&\Psi_{\lambda_{l_1}}^{\lambda_{l_3} V^{(l)}}(u)x=
\Delta(x)\Psi_{\lambda_{l_1}}^{\lambda_{l_3} V^{(l)}}(u),\qquad {\rm for}\  \forall x\in \dyslt.
\ena

We normalize the vertex operators as follows. 
\bea
&&\Phi_{\lambda_{l_1}}^{\lambda_{l_3} V^{(l)}}(u)
\ket{\lambda_{l_1}}=\ket{\lambda_{l_3}}\otimes w_m+
\cdots ,\label{normalI}\\
&&\Psi_{\lambda_{l_1}}^{\lambda_{l_3} V^{(l)}}(u)
\ket{\lambda_{l_1}}= w_m\otimes\ket{\lambda_{l_3}}+
\cdots \label{normalII}
\ena
with $m=(l+l_3-l_1)/2$. 
Expanding
$\Phi_{\lambda_{l_1}}^{\lambda_{l_3} V^{(l)}}(u)$ and 
$\Psi_{\lambda_{l_1}}^{\lambda_{l_3} V^{(l)}}(u)$ in a formal series
\bea
&&\Phi_{\lambda_{l_1}}^{\lambda_{l_3} V^{(l)}}(u)
=\sum_{m=0}^l\Phi_{l,m}(u)\otimes w_m,\label{expvertex}\\
&&\Psi_{\lambda_{l_1}}^{\lambda_{l_3} V^{(l)}}(u)
=\sum_{m=0}^l w_m\otimes\Psi_{l,m}(u),\label{expvertexII}
\ena
we obtain from (\ref{comulti}) and (\ref{evaluation}) the following intertwining relations. 

\vspace{5mm}
\noindent
{\bf Lemma 5.1}
\bea
{\rm Type\ I:}&&(l-m+1)\Phi_{l,m-1}(u)=[\Phi_{l,m}(u), f_0],\nonumber\\
&&[\Phi_{l,l}(u),e(v)]=0,\nonumber\\
&&h^+(v)\Phi_{l,l}(u)=
\frac{v-u-\frac{k-l+1}{2}\hbar}{v-u-\frac{k+l+1}{2}\hbar}
\Phi_{l,l}(u)h^+(v),\nonumber\\
&&h^-(v)\Phi_{l,l}(u)=
\frac{v-u+\frac{k+l-1}{2}\hbar}{v-u+\frac{k-l-1}{2}\hbar}
\Phi_{l,l}(u)h^-(v).\nonumber\\
{\rm Type\ II:}&&(m+1)\Psi_{l,m+1}(u)=[\Psi_{l,m}(u), e_0],\nonumber\\
&&[\Psi_{l,0}(u),f(v)]=0,\nonumber\\
&&h^{\pm}(v)\Psi_{l,0}(u)=
\frac{v-u-\frac{l-1}{2}\hbar}{v-u+\frac{l+1}{2}\hbar}
\Psi_{l,0}(u)h^{\pm}(v),\nonumber
\ena

\vspace{5mm}
\noindent
{\bf Lemma 5.2}
{\it The following vertex operators satisfy the intertwining relations in Lemma 5.1 in the limit $J\to \infty$ and $\delta\to 0$.}
\bea
&&\Phi_{l,m}(u)_{[J,\delta]}=\frac{1}{(l-m)!}
[\cdots[\Phi_{l,l}(u)_{[J,\delta]},
\underbrace{f_0]\cdots f_0]}_{l-m},\\
&&\Psi_{l,m}(u)_{[J,\delta]}=\frac{1}{m!}[\cdots[\Psi_{l,0}(u)_{[J,\delta]},
\underbrace{e_0]\cdots e_0]}_{m},
\quad (m=0,1,..,l)\ena
{\it with}
\bea
&&\Phi_{l,l}(u)_{[J,\delta]}\nonumber \\
&& =\exp\Bigl\{\sum_{j=0}^J[
\Phi^{(-)}(u;-\frac{k+l-1}{2}+2j;-\frac{k-l-1}{2}+2j)
+\Phi^{(-)}(0;2+\delta+2j;2j)]\Bigr\}\ 
\nonumber \\
&&\qquad \times e^{\frac{l}{k+2}a_{\Phi}}
\nonumber \\
&&\qquad\times \exp\Bigl\{-\sum_{j=1}^J
\Phi^{(+)}(u;-\frac{k-l+1}{2}-
(k+2)j;-\frac{k+l+1}{2}-(k+2)j)\Bigr\},
\\
&&\Psi_{l,0}(u)_{[J,\delta]}\nonumber \\
&& =\exp\Bigl\{\sum_{j=0}^J[
\Phi^{(-)}(u;\frac{l-5}{2}-k+2j;-\frac{l+1}{2}-k+2j)
+\Phi^{(-)}(0;2+\delta+2j;2j)]\Bigr\}\ 
\nonumber \\
&&\qquad \times e^{\frac{l}{k+2}a_{\Phi}}
\nonumber \\
&&\qquad\times \exp\Bigl\{-\sum_{j=1}^J
\Phi^{(+)}(u;-\frac{l-1}{2}-
(k+2)j;\frac{l-3}{2}-(k+2)j)\Bigr\}\nonumber\\
&&\qquad\times\exp{\Bigl\{\phi(u;\frac{l-5}{2}-k,-\frac{l+3}{2}-k)
+\chi(u;\frac{l-3}{2}-k)\Bigr\}},
\ena
where
\bea
&&\Phi^{(-)}(u;A;B)=\sum_{n>0}\frac{2a_{-n}}{(k+2)n}[(u+A\hbar)^n-
(u+B\hbar)^n]\\
&&\Phi^{(+)}(u;A;B)=
\sum_{n>0}\frac{a_{n}}{n}[(u+A\hbar)^n-
(u+B\hbar)^n]-\Bigl[\log\Bigl(\frac{u+A\hbar}{u+B\hbar}\Bigr)\Bigr]
\partial_{\Phi}\nonumber \\
\ena

\vspace{5mm}
\noindent
From this and the remark below the proposition 4.3, we arrive at the
following theorem.

\vspace{5mm}
\noindent
{\bf Theorem5.1}
{\it The vertex operators }
\bea
&&\Phi_{\lambda_{l_1}}^{\lambda_{l_3} V^{(l)}}(u)=\sum_{m=0}^l
\Phi^{(r)}_{l,m}(u)\otimes w_m,\\
&&\Psi_{\lambda_{l_1}}^{\lambda_{l_3} V^{(l)}}(u)=\sum_{m=0}^l
w_m\otimes \Psi^{(r)}_{l,m}(u),
\ena
{\it with}
\bea
&&\Phi^{(r)}_{l,m}(u)=g^{(I)\lambda_{l_3}V^{(l)}}_{\lambda_{l_1}}(u)
\oint_{{\cal C}}\frac{dt_1}{2\pi i}\oint\frac{dt_2}{2\pi i}
\cdots\oint\frac{dt_r}{2\pi i}
\Phi_{l,m}(u)_{[J,\delta]}S(t_1)_{[J]}
S(t_2)_{[J]}\cdots S(t_r)_{[J]},\nonumber\\
&&\Psi^{(r)}_{l,m}(u)=g^{(II)\lambda_{l_3}V^{(l)}}_{\lambda_{l_1}}(u)
\oint_{{\cal C}}\frac{dt_1}{2\pi i}\oint\frac{dt_2}{2\pi i}
\cdots\oint\frac{dt_r}{2\pi i}
\Psi_{l,m}(u)_{[J,\delta]}S(t_1)_{[J]}S(t_2)_{[J]}\cdots S(t_r)_{[J]},
\nonumber
\ena
{\it where $2r=l_1+l-l_3$ and 
$g^{(I,II)\lambda_{l_3}V^{(l)}}_{\lambda_{l_1}}(u)$ being the
 normalization function 
determined by the condition (\ref{normalI}) and (\ref{normalII}),
 gives the intertwiner 
in Definition5.1
 in the limit $J\to \infty$ and $\delta\to 0$.
The contour ${\cal C}$ is depicted in Fig.1\cite{Felder}.}

\section{Discussion}
We have constructed a free field representation of the level-$k$
Drinfeld currents for the Yangian double $\dyslt$ and screening operators.
As a result, we have obtained a deformation of the level-$k$ Wakimoto
modules. We also have realized the type I and type II vertex operators
intertwining these modules. 

A possible application of the results is a calculation of
correlation functions in massive integrable quantum field theory such
as higher spin $SU(2)$ invariant Thirring model and in higher spin XXX
spin chains. For this purpose, one has to make a precise identification
of the space of states with the deformed Wakimoto modules. It is also 
an important problem to derive the deformation of the
Knizhnik-Zamolodchikov equation\cite{Smir,BL} in the framework of free
field representation.

We hope to discuss these problems in future publication. 
 
\section{Acknowledgments} 
The author would like to thank Kenji Iohara for stimulating discussion.
He is also grateful to  Michio Jimbo,
Vladimir Korepin, Tetsuji Miwa and Max Niedermaier for discussions.
This work is supported by the COE research institute fellowship
and the Yukawa Memorial Foundation. 



\begin{thebibliography}{}
\bibliographystyle{unsrt}

\bibitem{Dr1}
Drinfeld, V.G. Hopf algebras and quantum Yang-Baxter equation.
{\it Soviet Math. Dokl.}{\bf 32} (1985), 254-258.

\bibitem{LS}
Leclair, A. and Smirnov, F.A. Infinite quantum group symmetry of
fields in massive 2D quantum field theory
{\it Int.J.Mod.Phys.}{\bf A7}(1992), 2997-3022.
 
\bibitem{Smir}
Smirnov, F.A. Dynamical symmetries of massive integrable models I, II
{\it Int.J.Mod.Phys.}{\bf A7} Suppl. {\bf 1B} (1992), 813-838, 839-858.

\bibitem{BL}
Bernard, D. and Leclair, A. The quantum double in integrable quantum field
theory
{\it Nucl.Phys.}{\bf B399}(1993), 709-748.

\bibitem{Dr2}
Drinfeld, V.G. Quantum groups, Proc. ICM Berkeley (1986), 798-820.


\bibitem{FR}
Frenkel, I.B. and Reshetikhin, N.Yu.
Quantum affine  algebras and holonomic difference equations.
{\it Commun. Math. Phys.}{\bf 146} (1992), 1-60.


\bibitem{JM}
Jimbo, M. and Miwa, T.
 Algebraic Analysis of Solvable Lattice Models.
 {\it Conference Board of the Math. Sci., Regional Conference Series
 in Mathematics},{\bf  85} (1995) and references therein.

\bibitem{IK}
Iohara, K. and Kohno, M. A central extension of ${\cal D}Y_{\hbar}(\glt)$
 and its vertex representations. UTMS 95-32, q-alg/9603032. 
to appear in {\it Lett. Math. Phys.}.

\bibitem{Ioh}
Iohara, K. Bosonic Representations of Yangian Double ${\cal D}Y_{\hbar}(g)$
with $g=gl_n, \sln$. Kyoto-Math 96-06, q-alg/9603033, to appear in
{\it Jour. of Phys. A}.

\bibitem{Kho}
Khoroshkin, S. Central Extension of the Yangian Double. 
q-alg/9602031.

\bibitem{KLP}
Khoroshkin, S., Lebedev, D., Pakuliak, S. Intertwining Operators for
the Central Extension of the Yangian Double. Preprint DFTUZ/95-28,
ITEP-TH-15/95, q-alg/9602030.

\bibitem{Dr3}
Drinfeld, V.G. A new realization of Yangians and quantized
affine algebras. {\it Soviet Math. Dokl.}{\bf  32} (1988), 212-216.

\bibitem{Wak}
Wakimoto, M.  Fock space representations of the affine Lie algebra 
$A^{(1)}_1$
{\it Comm. Math. Phys. }{\bf 104} (1986), 605-609.

\bibitem{Iohara}
Iohara, K. private communication.





\bibitem{Shir}
Shiraishi, J. Free boson representation of $U_q(\widehat{sl_2})$,
{\it Phys. Lett.} {\bf A171}(1992), 243-248.

\bibitem{Matsuo1}
Matsuo, A. Free field representation of quantum affine algebra 
$U_q(\widehat{sl_2})$, {\it Phys. Lett.} {\bf B308} (1992), 260-.

\bibitem{ABG}
Abada, A., Bougourzi, A.H. and Gradechi, El. Deformation of the Wakimoto construction, {\it Mod. Phys. Lett.} {\bf A8} (1993), 715-724.

\bibitem{KQS}
Kato, A., Quano, Y. and Shiraishi, J.
Free boson representation of $q-$vertex operators and their correlation functions,
{\it Comm. Math. Phys.} {\bf 157}(1993), 119-137.

\bibitem{Matsuo2}
Matsuo, A. A $q-$deformation of Wakimoto modules, primary fields and 
screening operators, {\it Comm. Math. Phys.} {\bf 160} (1994), 33-48.

\bibitem{Kon1}
Konno, H. BRST cohomology in quantum affine algebra 
$U_q(\widehat{sl_2})$,
{\it Mod. Phys. Lett.}{\bf A9}(1994),1253-1265.

\bibitem{Kon2}
Konno, H. Free field representation of quantum affine algebra 
$U_q(\widehat{sl_2})$ and form factors in higher spin XXZ model,
{\it Nucl. Phys.}{\bf B432}(1994), 457-486.

\bibitem{CP}
Chari, V. and  Pressley, A.
Yangians and $R-$matrices,  {\it L'Enseign.Math.t } {\bf 36} (1990), 
267-302. 
 
\bibitem{Felder}
Felder, G. BRST approach to minimal model,
{\it Nucl. Phys.}{\bf B317}(1989), 215-236.

\end{thebibliography}
\end{document}